\documentclass[12pt,letterpaper]{article}

\def\crseNo{}

\author{Michael La\,Croix
  \and
  Tom Roby
}
\title{Foatic actions of the symmetric group and fixed-point homomesy}
\usepackage{geometry}
\def\today{\number\day\space\ifcase\month\or
  January\or February\or March\or April\or May\or June\or
  July\or August\or September\or October\or November\or December\fi
  \space\number\year}
\parskip=\medskipamount
\usepackage{microtype}

\input{glyphtounicode} \pdfgentounicode=1

\usepackage{amsfonts}
\usepackage{amsmath}
\usepackage{amssymb}
\usepackage{bbm}
\newcommand{\II}{\mathbbm 1}
\usepackage{amsthm}
\usepackage[dvipsnames]{xcolor}
\newcommand{\cred}[1]{{\color{red}#1}}

\usepackage{graphicx}
\usepackage{fancyhdr}
\usepackage{listings}
\usepackage{manfnt}

\usepackage[normalem]{ulem}

\usepackage{url}

\usepackage{tikz-qtree}

\usepackage[colorlinks=true,citecolor=black,linkcolor=black,urlcolor=blue,
 pdfauthor={Tom Roby; Michael La Croix}]{hyperref}

\definecolor{lblue}{rgb}{0 .6 1}    
\usepackage{todonotes}    
\definecolor{bananamania}{rgb}{0.98, 0.91, 0.71}
\definecolor{apricot}{rgb}{0.98, 0.81, 0.69}
\def\Fix{\mathop{\rm Fix}}
\def\Rasc{\mathop{\rm Rasc}}
\def\exc{\mathop{\rm exc}}
\def\wexc{\mathop{\rm wexc}}
 
\def\multiset#1#2{\ensuremath{\left(\kern-.3em\left(\genfrac{}{}{0pt}{}{#1}{#2}\right)\kern-.3em\right)}}

\def\lcm{\mathop{\rm LCM}}

\def\eset{\emptyset}
\def\bull{\noindent $\bullet$\kern 2em}

\def\maj{\mathop{\rm{maj}}}

\def\bc#1#2{\left(\kern -2pt{#1\atop #2} \kern -2pt\right)}

\renewcommand{\hat}{\widehat}

\newcommand{\id}{\operatorname{id}}

\theoremstyle{plain}
\newtheorem{thm}{Theorem}
\newtheorem{lemma}[thm]{Lemma}
\newtheorem{prop}[thm]{Proposition}
\newtheorem{conj}[thm]{Conjecture}

\theoremstyle{definition}
\newtheorem{eg}[thm]{Example}
\newtheorem{defn}[thm]{Definition}

\newtheorem{prob.}[thm]{Problem}

\theoremstyle{remark}

\catcode`@=11
\def\m@th{\mathsurround\z@}
\def\cases#1{\left\{\,\vcenter{\normalbaselines\m@th
    \ialign{$##\hfil$&\quad##\hfil\crcr#1\crcr}}\right.}
\def\matrix#1{\null\,\vcenter{\normalbaselines\m@th
    \ialign{\hfil$##$\hfil&&\quad\hfil$##$\hfil\crcr
      \mathstrut\crcr\noalign{\kern-\baselineskip}
      #1\crcr\mathstrut\crcr\noalign{\kern-\baselineskip}}}\,}
\def\pmatrix#1{\left(\matrix{#1}\right)}
\def\bmatrix#1{\left[\matrix{#1}\right]}
\def\hang{\hangindent 24pt}
\def\d@nger{\medbreak\begingroup\clubpenalty=10000
  \def\par{\endgraf\endgroup\medbreak} %
  \noindent\hang\hangafter=-2
  \hbox to0pt{\hskip-\hangindent\dbend\hfill}}
\outer\def\danger{\d@nger}
\catcode`@=12 % at signs are no longer letters
\newcommand{\KK}{\mathbb K}

\newcommand{\calA}{\mathcal{A}}
\newcommand{\calB}{\mathcal{B}}

\newcommand{\cS}{{S}}
\newcommand{\cO}{\mathcal{O}}

\newcommand{\frakS}{\mathfrak{S}}

\def\inv{\mathop{\mathcal{I}}}
\def\comp{\mathop{\mathcal{C}}}
\def\divert{\mathop{\mathcal{D}}}
\def\rev{\mathop{\mathcal{R}}}
\newcommand{\rot}[1][2]{\mathop{\mathcal{Q}^{#1}}}
\def\qrot{\mathop{\mathcal{Q}}}
\def\foata{\mathop{\mathcal{F}}}
\def\fs{\mathop{\mathcal{S}}}
\def\opA{\mathop{\mathcal{A}}}
\def\opB{\mathop{\mathcal{B}}}
\begin{document}
%%This gets rid of the headers AND the rule in between the header and
 %%the live text. 
\lhead{}
\chead{}
\rhead{}
\renewcommand{\headrulewidth}{0pt}
\lfoot{\it \crseNo} 
\cfoot{\thepage}
\rfoot{\it \today }
\renewcommand{\footrulewidth}{0.4pt}
%\renewcommand{\marginparwidth}{0pt}
%\setlength{\marginparwidth}{1cm}

%%Put Class, Title, and Date below. 
%\hbox to \hsize{\large \textbf{Math  (Roby)} \hfil \textbf{XXXX} \hfil
%\bf Ongoing}   

\maketitle

\begin{abstract}

We study  maps on the set of permutations of $n$ generated by the R\'enyi-Foata map
intertwined with other dihedral symmetries (of a permutation considered as a 0-1 matrix).
Iterating these maps leads to dynamical systems that in some cases exhibit interesting
orbit structures, e.g., every orbit size being a power of two, and homomesic statistics
(ones which have the same average over each orbit).  In particular, the number of fixed
points (aka 1-cycles) of a permutation appears to be homomesic with respect to three of
these maps, even in one case where the orbit structures are far from nice.  For the most
interesting such ``Foatic'' action, we give a heap analysis and recursive structure that
allows us to prove the fixed-point homomesy and orbit properties, but two other cases
remain conjectural.
\end{abstract}

%% \tableofcontents

%\setcounter{section}{-1}
\section{Introduction}\label{sec:int}

\subsection{The R\'enyi-Foata Map on permutations}\label{ss:RFmap}
% *** Discuss background on the R\'enyi-Foata map $\foata $ that drops parentheses, and give quick
% definition.  

A well-known bijection $\foata$ from the symmetric group $\frakS_{n}$ to itself, due
to R\'enyi~\cite[\S 4]{Ren62} and Foata-Sch\"utzenberger~\cite[pp. 13--15]{FS70}, simply takes
a permutation given in a canonical disjoint cycle decomposition, drops the parentheses, and
reinterprets the result as a permutation in one-line notation.  (See (\ref{eq:foata}) for an example.)  
Although not respecting the algebraic structure of $\frakS_{n}$, it provides insight into
combinatorial properties of permutations.  For example it shows that among
permutations in $\frakS_{n}$, the number with exactly $k$ \emph{cycles} (counted by signless Stirling numbers
of the first kind) is the same as the number with exactly $k$ \emph{left-to-right maxima},
which we call here \textbf{records} for short,  and that the
number with $k$ \emph{ascents} (counted by Eulerian numbers) is the same as the number with
$k+1$ \emph{weak excedances}.  

\begin{defn}\label{def:RFmap}
Let $w \in \frakS_{n}$.  The \textbf{canonical (disjoint) cycle decomposition (CCD)} of
$w$ is the decomposition of the bijection $w$ into disjoint cycles, where (a) each cycle
is written with its largest element first and (b) the cycles are written in increasing order
of first (largest) elements.  The \textbf{R\'enyi--Foata map} $\foata: \frakS_{n}\rightarrow \frakS_{n}$ simply
removes the parentheses from the CCD of $w$ and regards the resulting word as a permutation in
one-line notation.  
\end{defn}

For example: 
\begin{equation}\label{eq:foata}
w = 847296513 = (4 2)(6)(8 1)(9 3 7 5) \stackrel{\foata}{\mapsto} 426819375 = (2)(9 5 1 4 8 7 3 6 ),
\end{equation}

where the input and output to $\foata$ are each written both in one-line notation and CCD.  
Note that here $w$ has 4 cycles, and $\foata (w)$ has 4 records (viz., 4, 6, 8,
and 9).  

It is easy to see that $\foata$ is a bijection, whose inverse is given as follows. Take a
permutation $w$ given in one-line notation, place a left parenthesis before each record, then
place corresponding right parentheses (one before each internal left parenthesis
and after the last element).  

Every permutation $w$ has multiple representations, and we use whichever one is
convenient at the moment.  Ironically, from a computer science (or even just hand
computational) standpoint, the conversion between different representations of the same $w$
requires more effort than the map $\foata$ itself. 

For a broader overview of this corner of permutation enumeration and its history, we direct the reader to the
text of Stanley~\cite[\S 1.3]{ec1ed2} and the discussion
of references in that chapter. Therein our map
$w \mapsto \foata(w)$ is called the \textit{fundamental bijection} (translating the
term ``transformation fondamentale'' of \cite{FS70}) and denoted $w\mapsto \hat{w}$. Bona's
text~\cite[\S3.3.1]{BonaPerm} is another useful reference.

\subsection{Dynamics of permutation bijections}\label{ss:dyn}

We consider certain cyclic actions on the symmetric group that are generated as an intertwining of
the R\'enyi-Foata map $\foata$ with other involutive dihedral symmetries on $\frakS_{n}$
(via permutation matrices). In particular, we consider actions generated by
maps of the following form: 
\begin{equation}\label{eq:foatic}
\frakS_{n} 
\stackrel{\foata}{\rightarrow} 
\frakS_{n} 
\stackrel{\calA}{\rightarrow}
\frakS_{n} 
\stackrel{\foata^{-1}}{\rightarrow}
\frakS_{n} 
\stackrel{\calB}{\rightarrow} 
\frakS_{n} 
\end{equation}
where $\calA$ and $\calB$ are dihedral involutions, defined below.  

\begin{defn}\label{def:invs}
For completeness, we name the five \textbf{dihedral involutions} of
$\frakS_{n}$ as follows: 
\begin{enumerate}
\def\theenumi{\alph{enumi}}\def\labelenumi{(\theenumi)}
\item $\comp: \frakS_{n}\rightarrow \frakS_{n}$, which takes a permutation $w = w_{1}\dots w_{n}$ to its \textbf{complement}
whose value in position $i$ is $n+1-w_{i}$; 

\item $\rev: \frakS_{n}\rightarrow \frakS_{n}$, which takes a permutation $w = w_{1}\dots w_{n}$ to its \textbf{reversal}
whose value in position $i$ is $w_{n+1-i}$; 

\item $\rot : \frakS_{n}\rightarrow \frakS_{n}$, which takes a permutation $w = w_{1}\dots w_{n}$
to its \textbf{rotation by 180-degrees}, whose value in position $i$ is $n+1-w_{n+1-i}$.  (We
reserve $\qrot$ to denote the dihedral symmetry \textbf{rotation of (the permutation matrix
of) $w$ by 90-degrees counterclockwise}.)   

\item $\inv: \frakS_{n}\rightarrow \frakS_{n}$, which takes a permutation $w$ to its \textbf{inverse} $w^{-1}$; 

\item $\divert: \frakS_{n}\rightarrow \frakS_{n}$, which takes a permutation $w$ to its
its \textbf{rotatedInverse} $\rot (\inv (w))$. 
\end{enumerate}
\end{defn}

\begin{figure}
  \caption{\label{fig:DihedralActions}The dihedral
      symmetries on $\frakS_n$} 
\medskip
\mbox{}\hspace{\fill}
\begin{minipage}{.45\textwidth}
\includegraphics[width=\textwidth]{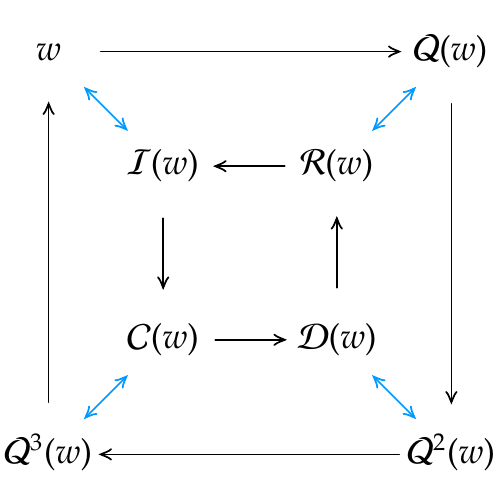}%
\end{minipage}
\hspace{\fill}
\begin{minipage}{.45\textwidth}
\includegraphics[width=\textwidth]{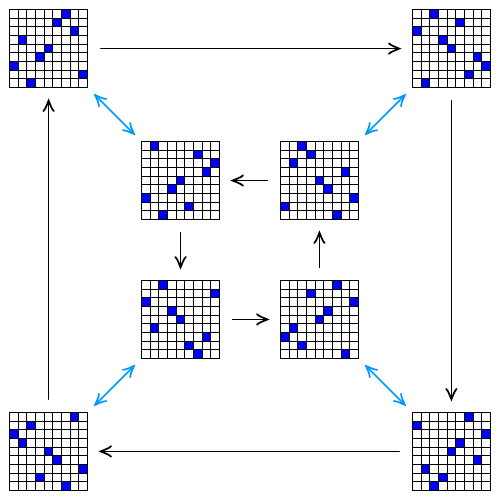}%
\end{minipage}
\hspace{\fill}\mbox{}\\[1.8ex]
\textit{\textcolor{black}{The left image shows the dihedral group
    presented in terms of generators $\rot[]$ (black) and
    \textcolor{lblue}{$\inv$ (blue)}.  Its actions the permutation
    graph of $w=361458972$ are shown on the right.  Note that $\rot[]$
    corresponds to a counter-clockwise rotation of permutation
    matrices, but this induces a clockwise rotation of the
    corresponding graphs.}}
\end{figure}

Of course, some of these maps can be obtained as compositions of others, e.g., $\rot = \comp
\circ \rev = \rev \circ \comp$.   \textcolor{black}{Figure~\ref{fig:DihedralActions}
summarizes the relationships between the dihedral involutions by
presenting each in terms of the generators $\rot[]$ (black) and
\textcolor{lblue}{$\inv$ (blue)}.}

For humans, the first two operations ($\comp$ and $\rev$) are easy to
compute on permutations in one-line notation, but not in CCD.  Taking
the inverse of $w$ in CCD is easy, since one just reverses the order
of non-maximal elements in each cycle, but harder in one-line
notation.  $\rot (w)$ is easy to compute in CCD: one simply
complements each element of each cycle.  The result will be a
cycle-decomposition that is not CCD in general, but is easy to convert
to CCD.  In one-line notation $\rot (w)$ also easy: one reads off the
complement of each element of $w$ from right-to-left.
\textcolor{black}{The operators $\rev$, $\comp$, and $\rot$ act naturally on one-line
  presentation of permutations, while $\inv$, $\divert$, and $\rot$
  act naturally on the CCD.}

\begin{defn}\label{def:Foatic}
We call the invertible maps defined as the fourfold composition in
(\ref{eq:foatic}), where $\calA$ and $\calB$ are from the above list,
\textbf{Foatic}.  This gives a total of 25 different Foatic maps to
study.
\end{defn}

\textcolor{black}{As part of our investigation we occasionally permitted
  $\calA$ and $\calB$ to act as $\rot[]$, $\rot[3]$, or the identity,
  but these additional actions did not produce any interesting results
  from the perspective of our investigation into homomesy.  By
  excluding $\rot[]$ and $\rot[3]$ from our definition, we preserve
  the potentially useful fact that every Foatic action factors as a
  product of two involutions.  In fact, the most tractable actions
  seemed to occur only when $\calA$ acted naturally on one-line
  presentations and $\calB$ acted naturally on CCD, although Sheridan
  Rossi later identified statistics of interest for all five
  possibilities of $\calA$ and $\calB$ summarized in Appendix~A of
  \cite{esrthesis}.}

\begin{eg}\label{eg:invs}

Let $w=361458972$ (in one-line notation) $=(31)(4)(5)(92687)$ in CCD.  Then

\begin{minipage}{.6\textwidth}
\begin{enumerate}
\def\theenumi{\alph{enumi}}
\def\labelenumi{(\theenumi)}

\item $\comp (w) = 749652138 = (5)(624)(71)(983)$, 
      
\item $\rev (w) = 279854163 = (5)(648)(712)(93)$, 
      
\item $\rot (w) = 831256947 = (5)(6)(84231)(97)$,  

\item $\inv (w) = 391452867 = (31)(4)(5)(97862)$, and

\item $\divert (w) = 342856917= (5)(6)(81324)(97)$.

\end{enumerate}
\end{minipage}%
\begin{minipage}{0.38\textwidth}
  \includegraphics{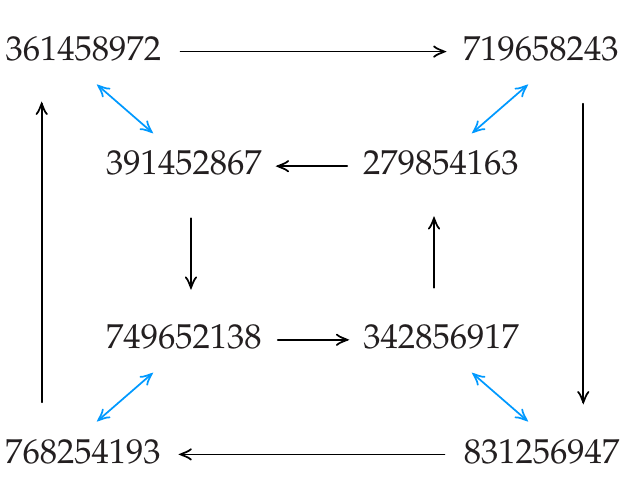}
\end{minipage}

\end{eg}

It is worth noting that any Foatic map can be thought of as a product of two involutions
$\opB \circ \foata^{-1}\circ \opA \circ \foata = \opB (\foata^{-1}
\opA  \foata )$, where one of our dihedral involutions has been conjugated by
$\foata$.  It will also be convenient in some cases to consider these actions to be starting
partway through the composition, e.g., as $\foata^{-1} \opA  \foata
\opB$\textcolor{black}{, considering such a conjugate action will produce the same orbit
structure, but with potentially different homomesic statistics}.

\begin{eg}\label{eg:foatic}
If $\calA =\comp $ and $\calB =\inv $, then we get the Foatic map $\gamma  : = \inv \circ
\foata^{-1}\circ \comp \circ \foata$.  If $n=5$, then $\gamma  [(4213)(5)] = (2)(4)(513)$ as follows
\[
w = (4213)(5)\stackrel{\foata}{\mapsto} 42135 
	 \stackrel{\comp}{\mapsto}  24531
	 \stackrel{\foata^{-1}}{\mapsto} (2)(4)(531)
	 \stackrel{\inv }{\mapsto} (2)(4)(513) = \gamma (w).
\]
The orbit (of size six) generated by the above $w$ is
\[
w = (4213)(5)\stackrel{\gamma }{\mapsto} (2)(4)(513) 
	 \stackrel{\gamma }{\mapsto} (412)(53)
	 \stackrel{\gamma }{\mapsto} (2)(5314)
	 \stackrel{\gamma }{\mapsto} (431)(52)
	 \stackrel{\gamma }{\mapsto} (2)(3)(541)
	 \stackrel{\gamma }{\mapsto}  w
\]

\end{eg}

\subsection{The homomesy phenomenon}\label{ss:homomesy}

The homomesy phenomenon was first isolated by Propp and the second author~\cite{propp-roby}
around 2011.  Given a group action on a set of 
combinatorial objects, a statistic on these objects is called \textbf{homomesic} if its
average value is the same over all orbits. More precisely: 

\begin{defn}\label{defn:homomesic}
Given a set $\cS$, %%% of combinatorial objects, 
an invertible map $\tau$ from $\cS$ to itself such that each
$\tau$-orbit is finite, 
and a function (or ``statistic'') $f: \cS \rightarrow \KK$
taking values in some field $\KK$ of characteristic zero,
we say the triple $(\cS,\tau,f)$ exhibits {\bf homomesy}
iff there exists a constant $c \in \KK$
such that for every $\tau$-orbit $\cO \subset \cS$
\begin{equation}
\label{general-ce}
\frac{1} {\#\cO}
\sum_{x \in \cO} f(x) = c .
\end{equation}
In this situation
we say that the function $f: \cS \rightarrow \KK$ is 
{\bf homomesic}
under the (cyclic) action of $\tau$ on $\cS$,
or more specifically \textbf{\textit{c}-mesic}.
\end{defn}

When $\cS$ is a finite set,
homomesy can be restated equivalently as 
all orbit-averages being equal to the global average:

There are many examples of this phenomenon,
of varying degrees of difficulty, in both older and more recent combinatorially work;
see~\cite{robydac} for a survey.  In particular, in many cases it has been fruitful to 
investigate actions that can be viewed as the composition of involutions on the set of
objects.  These include \emph{promotion} of semi-standard Young tableaux~\cite{BPS13} and \emph{rowmotion}
on the set of antichains or order ideals of certain posets~\cite{ast, strikerwilliams,
rowmotion-in-slow-motion}, (particularly \emph{minuscule})
posets~\cite{rush-wang-homomesy-minuscule}. 
Extensions include
looking at certain more general products of \emph{toggling} involutions (see ~\cite{strikergentog}) in a variety of
contexts, including non-crossing partitions~\cite{aim-nctoggles}, independent sets of a path
graph~\cite{tispg}, antichains~\cite{joseph-antichain-toggling}.  Liftings of these
combinatorial maps to the piecewise-linear setting of polytopes and further (via
detropicalization) to the birational (and even noncommutative) realms are also of interest 
~\cite{einpropp,GrRo1,GrRo2,BAR-motion}. 

While homomesic statistics can often be found in situations which give interesting examples
of the \emph{cyclic sieving
phenomenon} of Reiner, Stanton, and White~\cite{csp,csp-brief}, they also turn up for actions
whose orbit structure is not well behaved.  In particular, data strongly suggests that the
Foatic complement-rotation map $\rho$ of Section~\ref{sec:comp-rot} has homomesic fixed-point
statistic despite its orbit structure.  The ``Coxeter-toggling'' of independent sets of a
path graph~\cite{tispg} is another example with proven homomesy but where the orbit
structures are too messy to find nice cyclic sieving. 

\subsection{Summary of the paper}\label{ss:summ}

This paper had its origin in James Propp's vision to look for natural homomesies among basic
combinatorial objects, particularly those counted by Rota's Twelvefold
Way~\cite[\S1.9]{ec1ed2}.  In this first section we give background and define the basic
setup.  In Section~\ref{sec:rev-inv}, we give a careful study of reversal-inclusion, the
most interesting Foatic map.  In this situation, we use \emph{heaps}, a kind of decreasing binary
tree, to understand the recursive nature of the map.  This allows us to prove some
properties of the orbit structure, and find several homomesic statistics, including 
the natural statistic $\Fix $ on $\frakS_{n}$ which counts the number of fixed points
(1-cycles) of a permutation. 

In Sections~\ref{sec:comp-inv}, \ref{sec:comp-rot},
and~\ref{sec:rev-rot} we consider \textcolor{black}{three} other Foatic
maps, complement-inversion, complement-rotation, and
reversal-rotation, which appear conjecturally to have $\Fix$ as a
homomesic statistic.  The former also seems to have nice orbits.  We
give data to support the conjectures which remain open.  Finally, in
Section~\ref{sec:otherMaps}, we note that computer-generated data
shows than none of the other \textcolor{black}{21} possible Foatic maps
exhibits homomesy for $\Fix$.

In her doctoral dissertation (supervised at UConn by the second author), Elizabeth Sheridan
Rossi searched for homomesy for (linear combinations of) a wider range of permutation
statistics under Foatic actions.  For example, let $\Fix_{i}(w)$ be the indicator function
that takes the value 1 if $w_{i}=i$ and 0 otherwise.  Then $\Fix_{1}-\Fix_{n}$ is homomesic
for 13 of the 25 Foatic maps.  Call $i\in [n]$ an \textbf{excedance} if $w_{i}>i$ and a
\textbf{weak excedance} if $w_{i}\geq i$; set $\exc w =$ number of excedances and
$\wexc = $ number of weak excedances of $w\in \frakS_{n}$.  Sheridan Rossi proves that
$\wexc +\exc$ is homomesic for inversion-inversion, rotation-inversion and
rotation-rotation.  She conjectures that $\wexc$ is homomesic for complement-rotatedInverse
and reversal-rotatedInverse. See Appendix~A for a summary and Chapter~2 of~\cite{esrthesis}
for more information.  

In the same dissertation, Sheridan Rossi considered a variant on the above setup, replacing the
R\'enyi--Foata map $\foata $ with the Foata--Sch\"utzenberger map $\fs$.  The latter map is a
(somewhat complicated) bijection on $\frakS_{n}$ with the property that $\maj (w)= \inv
\fs(w)$.  It gave the first bijective proof that the statistics $\maj$ and $\inv$ are
equidistributed on $\frakS_{n}$ (a result first obtained by MacMahon).  See~\cite[\S
1.4]{ec1ed2} for further background. Let $D_{i}(w)$ be the indicator function of whether $i\in
[n]$ is a descent of $w\in \frakS_{n}$. Among other homomesies, Sheridan Rossi proved that
$D_{1}+D_{n-1}$ is homomesic for five such intertwinings, and conjectured that this extends
to $D_{i}+D_{n-i}$ for $1<i<n$.  See Appendix~B for a summary and Chapter~3 of~\cite{esrthesis}
for more information.

\subsection{Acknowledgments}\label{ss:ack}

The authors are particularly grateful to James Propp, who first suggested this line of
inquiry and noticed that the number of fixed points is a homomestic statistic for a few of the Foatic
maps we consider here.  We also greatly appreciate conversations we've had with Ira
Gessel, and Darij Grinberg.  David Einstein, Michael Joseph, and Elizabeth Sheridan Rossi
all read early drafts of this paper, finding errors and contributing insightful comments.
Exploratory computations were carried out in the PostScript programming language, and later
confirmed by Sage code written by David Einstein.  This collaboration began when the first
author was a postdoc and the second a research affiliate in the mathematics department at
MIT, whose hospitality we gratefully acknowledge. 

\section{Reversal-inversion}\label{sec:rev-inv}

The Foatic action with the nicest orbit structures and properties is the one given by:
\begin{equation}\label{eq:phibar}
\overline{\varphi}: \frakS_{n} 
\stackrel{\foata}{\rightarrow} 
\frakS_{n} 
\stackrel{\rev}{\rightarrow}
\frakS_{n} 
\stackrel{\foata^{-1}}{\rightarrow}
\frakS_{n} 
\stackrel{\inv}{\rightarrow} 
\frakS_{n}. 
\end{equation}
It turns out to be easier to study its conjugate map: 
\begin{equation}\label{eq:phi }
\varphi: \frakS_{n} 
\stackrel{\rev}{\rightarrow}
\frakS_{n} 
\stackrel{\foata^{-1}}{\rightarrow}
\frakS_{n} 
\stackrel{\inv}{\rightarrow} 
\frakS_{n} 
\stackrel{\foata}{\rightarrow} 
\frakS_{n} 
\end{equation}
and to use the representation of permutations as \emph{heaps}, aka \emph{decreasing binary
trees}~\cite[\S1.5]{ec1ed2}.  The action of $\varphi $ (equivalently $\overline{\varphi}$)
transforms the heap in an easily described way that preserves the isomorphism class of the
underlying unlabeled tree.  As a result we get a simple recursive structure of how $\varphi$
acts on permutations of $n$ in one-line notation, i.e., $\varphi (AnB) = Bn\varphi (A)$,
where $A$ and $B$ are partial permutations (also in one-line notation).  This allows us to
prove that the statistic $\Fix$ on $\frakS_{n}$ that counts the number of fixed points is
1-mesic, and also show that all the orbit sizes are powers of 2 (Theorem~\ref{thm:phi}).  

\begin{figure}
\begin{center}
  \caption{\label{fig:phiOrbits}Two orbits (one for $\frakS_{7}$, one
    for $\frakS_{9}$) of the Foatic reversal-inclusion map
    $\overline{\varphi}$ with associated heaps, with \cred{fixed
      points} marked in \cred{red}.  Each orbit has an average of one
    fixed point per permutation. }

\bigskip

\includegraphics[width=4.2in]{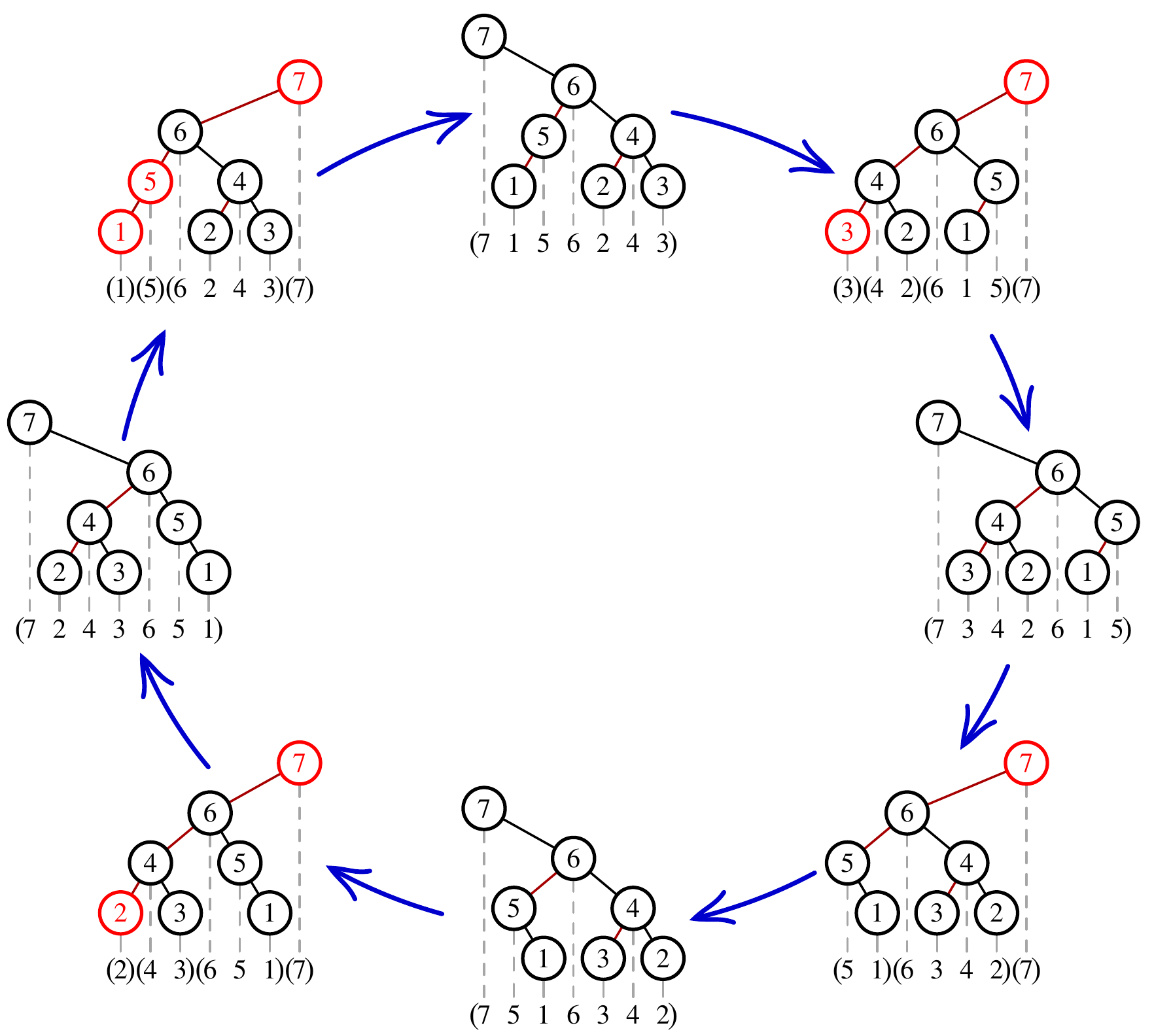}
\vskip 0.6in

\includegraphics[width=4.2in]{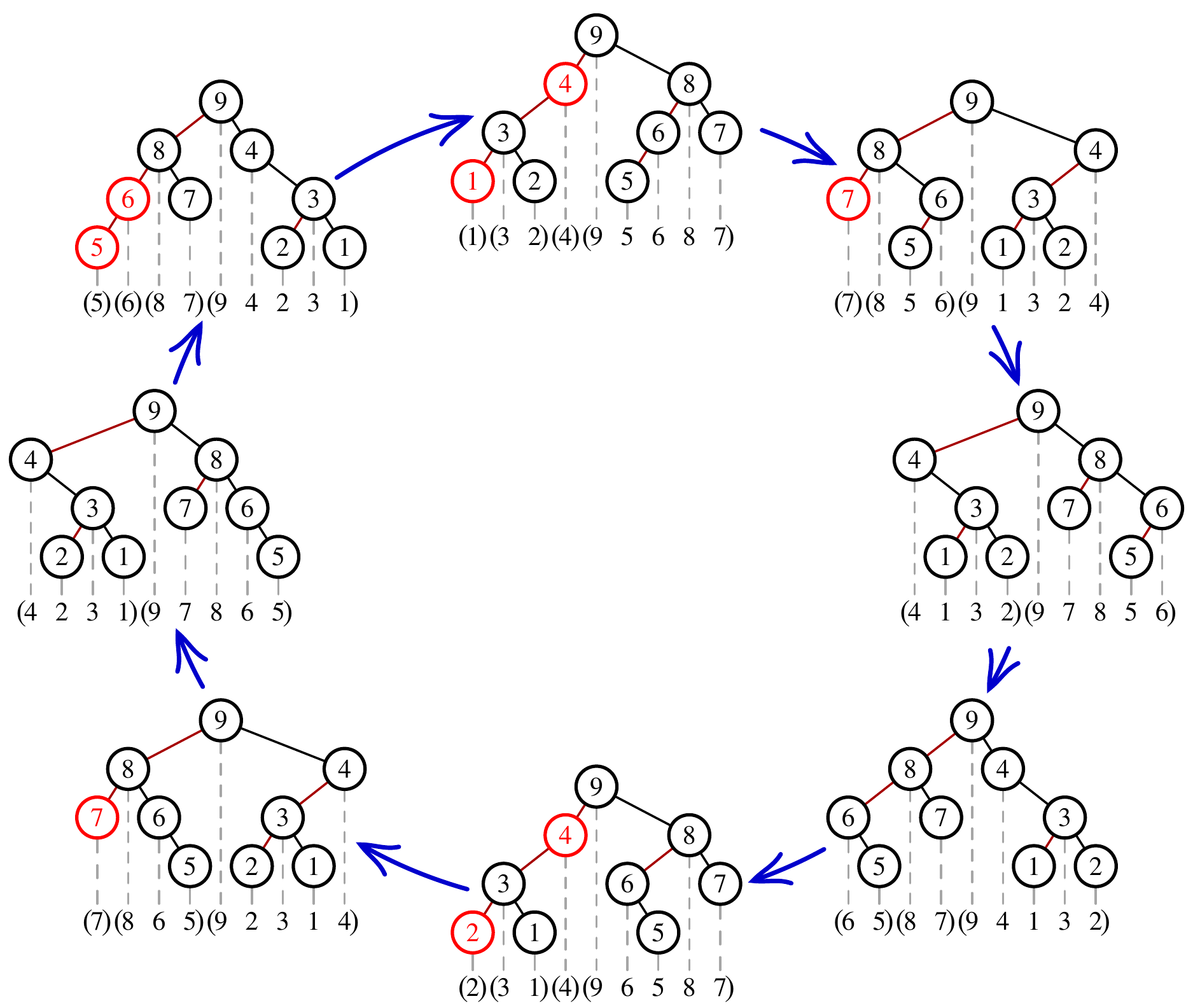}
\end{center}
\end{figure}

\begin{eg}\label{eg:phi} Let $w=(2)(43)(51)\in \frakS_{5}$ in CCD.  Then the successive action of
$\overline{\varphi}$ on $w$ is detailed below: 
%%% Note that alignat alternates between flushright and flushleft, hence some & need to be doubled
\begin{alignat*}{5}\label{eqn:phiEG}
w = &(2)(43)(51) && \mapsto 24351 \mapsto 15342 \mapsto (1)(5342)   &&\mapsto (1)(5243)   && = \overline{\varphi}(w)\\ 
    &(1)(5243)   && \mapsto 15243 \mapsto 34251 \mapsto (3)(42)(51) &&\mapsto (3)(42)(51) && = \overline{\varphi}^{2}(w)\\
    &(3)(42)(51) && \mapsto 34251 \mapsto 15243 \mapsto (1)(5243)   &&\mapsto (1)(5342)	 && = \overline{\varphi}^{3}(w)\\
    & (1)(5342)  && \mapsto 15342 \mapsto 24351 \mapsto (2)(43)(51) &&\mapsto (2)(43)(51) && = \overline{\varphi}^{4}(w)
\end{alignat*}
This example also displays (down the second column) the conjugate orbit of $\varphi$, also of size 4.  
\[
24351 \stackrel{\varphi}{\rightarrow}
15243 \stackrel{\varphi}{\rightarrow}
34251 \stackrel{\varphi}{\rightarrow}
15342 \Lsh 
\]
\end{eg}
Larger examples with two orbits from the $\overline{\varphi}$-actions on $\frakS_{7}$ and $\frakS_{9}$
are given in Figure~\ref{fig:phiOrbits}.  

% \begin{figure}
% \caption{An orbit of $\inv \foata^{-1}\rev \foata $}
% \label{fig:rev-inv orbit}
% Should we give it in one-line notation, then switch to cycles/heaps later? 
% \end{figure}

\begin{defn}\label{def:heap}
Let $S$ be a finite totally ordered set, and $w\in \frakS_{S}$ a permutation of $S$ written
in one-line notation.  If $S\subseteq [n]: = \{1,2,\dots ,n \}$, we call $w$ a
\textbf{partial permutation of $\mathbf{n}$}.  We recursively define the \textbf{heap} of $w$, $H(w)$ 
as follows.  Set $H(\eset \text{ (the empty word)}) = \eset$ (the empty tree). If $w\neq \eset$, let $m$ be the
largest element of $w$, so $w$ can be written uniquely as $umv$, where $u$ and $v$ are
partial permutations (possibly empty).  Set $m$ to be the root of $H(w)$, with $H(u)$ its
left subtree and $H(v)$ its right subtree.  
\end{defn}
The heap of a permutation will turn out to be a \emph{decreasing binary tree}, i.e., the labels along any
path from the root form a decreasing sequence.  For more information, see the equivalent
definition of \emph{increasing binary tree} and the results thereafter
in~\cite[\S1.5]{ec1ed2}.  

\begin{eg}\label{eg:heap}
The heap associated with $w=314975826$ is shown below. 
\begin{center}
\tikzset{every tree node/.style={minimum width=0.5em,draw,circle},
         blank/.style={draw=none},
         edge from parent/.style=
         {draw, edge from parent path={(\tikzparentnode) -- (\tikzchildnode)}},
         level distance=1.2cm,sibling distance=.5cm}
\begin{tikzpicture}[scale=0.8]
\Tree
[.9     
    [.4 
   	[.3 
           \edge[blank]; \node[blank]{}; 
           [.1 ]
        ] 
     \edge[blank]; \node[blank]{}; 
    ]	
    [.8  
       [.7 
       	   \edge[blank]; \node[blank]{}; 
           [.5 ]
       ] 
       [.6 
           [.2 ]
           \edge[blank]; \node[blank]{};
         ]
    ]
]
\end{tikzpicture}
\end{center}
 \end{eg}
% \begin{tikzpicture}[every tree node/.style={draw,circle},
%    level distance=1.25cm,sibling distance=.5cm, 
%    edge from parent path={(\tikzparentnode) -- (\tikzchildnode)}]
% \Tree [.\node[red] {9}; 
%     [.4  
%       [.4  ] [.1 ] [.3 ]
%     ]
%     [.8
%       [.6 ] 
%     ] ]
% \end{tikzpicture}
% \[
% \xymatrixcolsep{0.4pc}
% \xymatrixrowsep{0.6pc}
% \xymatrix{
% && \fullmoon \ar@{-}[ld] \ar@{-}[rd] & & \fullmoon \ar@{-}[ld] \ar@{-}[rd] & \\
% \rho_{J}: & \newmoon \ar@{-}[rd] & & \fullmoon \ar@{-}[ld] \ar@{-}[rd] & & \fullmoon \ar@{-}[ld]& \longrightarrow  \\
% && \newmoon & & \newmoon &
% }
% \] 

We state the following very elementary facts without proof. 

\begin{prop}
\label{prop:heaps}
Let $w\in \frakS_{n}$ given in one-line notation have corresponding heap $H(w)$. 
\begin{enumerate}

\item In $H(w)$, the left successor of $j$ is the greatest element $k$ to the left of $j$ in $w$, such
that all elements of $w$ between $k$ and $j$ inclusive are $\leq j$

\item The map $w\mapsto H(w)$ is a bijection between $\frakS_{n}$ and the set
of decreasing binary trees (as defined above) with $n$ vertices.  

\item Let $\sigma\in \frakS_{S}$, where $S=\{x_{1} < x_{2} <\dots < x_{\ell} \}$ is any
finite totally ordered set.  Let $\xi : S\rightarrow [\ell ]$ be the canonical bijection
$x_{i}\mapsto i$, which naturally extends to a bijection $\xi : \frakS_{S}\rightarrow
\frakS_{\ell}$.  We can use this map to extend the notions of dihedral symmetries, cycle
structure, and the R\'enyi-Foata map $\foata$ to $\frakS_{S}$, in particular to partial
permutations of $n$.  

\end{enumerate}

\end{prop}

We need the last statement above to be able to define and prove things recursively.  

\begin{eg}\label{eg:xi}
Consider the partial permutation $\sigma = 753296$ of 9.  We can consider this in two-line
notation as $\vphantom{\bmatrix{T^{w}\cr R\cr V}}\sigma = \pmatrix{2&3&5&6&7&9\cr 7&5&3&2&9&6}$, and in CCD as  $\sigma
=(53)(9627)$.  
Here $\sigma^{-1} = \pmatrix{2&3&5&6&7&9\cr 6&5&3&9&2&7} 
 = 654927 = (53)(9726)$.  
\end{eg}

\begin{lemma}\label{lem:recurPhi}
Let $w\in \frakS_{n}$ have the form $AnB$ (in one-line notation), where $A$ and $B$ are (possibly empty) partial
permutations of $n$.  Then the action of $\varphi$ satisfies $\varphi (AnB) = \varphi (B)n
A$.  Thus, $H(\varphi (AnB))$ is the heap interchanging the left and right subtrees at the
root vertex $n$, leaving
the former unchanged and applying $\varphi$ recursively to the latter.  
In particular, the action of $\varphi$
preserves the underlying unlabeled graph of the corresponding heaps.  
\end{lemma}
\begin{proof}
We analyze step-by-step how $\varphi = \foata\circ \inv \circ \foata^{-1} \circ \rev$ acts on the word
$w=AnB$, using $w=314975826$ (from Example~\ref{eg:heap}) as a running example.  First we see
straight from the definitions that $\rev$ acts on the heap by switching left and right at
each vertex, yielding the middle diagram in Figure~\ref{fig:heap}.   
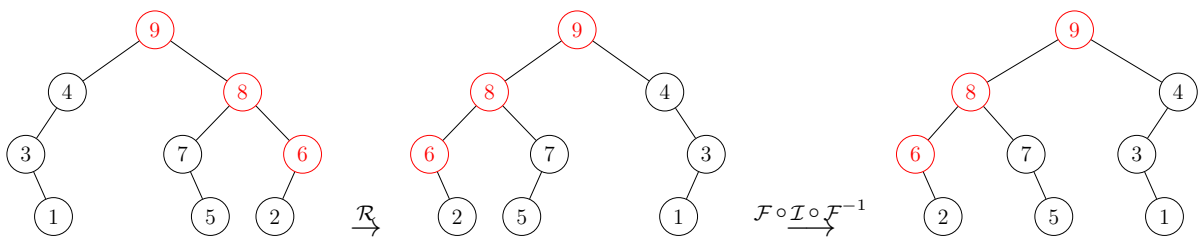
\begin{figure}[h]
\caption{The action of reversal-inclusion on the heap of the permutation $w=314975826$.}
\label{fig:heap}
\bigskip
\begin{center}
\tikzset{every tree node/.style={minimum width=0.5em,draw,circle},
         blank/.style={draw=none},
         edge from parent/.style=
         {draw, edge from parent path={(\tikzparentnode) -- (\tikzchildnode)}},
         level distance=1.2cm,sibling distance=.5cm}
\begin{tikzpicture}[scale=0.69]
\Tree
[.\node[red] {9};     
    [.4 
   	[.3 
           \edge[blank]; \node[blank]{}; 
           [.1 ]
        ] 
     \edge[blank]; \node[blank]{}; 
    ]	
    [.\node[red] {8};  
       [.7 
       	   \edge[blank]; \node[blank]{}; 
           [.5 ]
       ] 
       [.\node[red] {6}; 
           [.2 ]
           \edge[blank]; \node[blank]{};
         ]
    ]
]
\end{tikzpicture}
% \begin{tikzpicture}[scale=0.69]
% \Tree
% [.9     
%     [.4 
%    	[.3 
%            \edge[blank]; \node[blank]{}; [.1 ]
%         ] \edge[blank]; \node[blank]{}; 
%     ]	
%     [.8  
%        [.7 
%        	   \edge[blank]; \node[blank]{}; [.5 ]
%        ] 
%     \edge[]; [.6 
%              \edge[]; {2}
%              \edge[blank]; \node[blank]{};
%          ]
%     ]
% ]
% \end{tikzpicture}
$\stackrel{\rev}{\rightarrow}$
\begin{tikzpicture}[scale=0.69]
\Tree
[.\node[red] {9};    
    [.\node[red] {8};      
       [.\node[red] {6};    
             \edge[blank]; \node[blank]{};
             \edge[]; {2}
       ]
       [.7 
           [.5 ]
       	   \edge[blank]; \node[blank]{}; 
       ] 
    ]
    [.4 
	\edge[blank]; \node[blank]{}; 
   	[.3 
           [.1 ]
           \edge[blank]; \node[blank]{}; 
        ]  
    ]	
]
\end{tikzpicture}
$\stackrel{\foata \circ \inv \circ \foata^{-1} }{\longrightarrow}$
\begin{tikzpicture}[scale=0.69]
\Tree
[.\node[red] {9};    
    [.\node[red] {8};      
       [.\node[red] {6};    
             \edge[blank]; \node[blank]{};
             \edge[]; {2}
       ]
       [.7 
       	   \edge[blank]; \node[blank]{}; 
           [.5 ]
       ] 
    ]
    [.4 
   	[.3 
           \edge[blank]; \node[blank]{}; 
           [.1 ]
        ]  
	\edge[blank]; \node[blank]{}; 
    ]	
]
\end{tikzpicture}
\end{center}
\end{figure}

Here we mark in \cred{red} the sequence of records of $\rev (w) =
\cred{6}2\cred{8}57\cred{9}413$, which were the 
original right-to-left maxima in $w$.  These will be exactly those elements starting a cycle
(following a
left-parenthesis) in $\foata^{-1}(\rev (w)) = (\cred{6}2)(\cred{8}57)(\cred{9}413)$.  Applying
$\inv$ to this 
(maintaining CCD) reverses the elements to the right of the first/largest element within
each cycle, e.g., $\inv (\foata^{-1}(\rev (w))) = (\cred{6}2)(\cred{8}75)(\cred{9}314)$.  Finally, applying
$\foata$ drops parentheses, leaving the same set of records as before (for $\rev (w)$),
e.g., $\varphi (w) = \cred{6}2\cred{8}75\cred{9}314$.  

In the subtree of the original left subword $A=314$,
left and right were switched by $\rev$, then back again by $\foata \circ \inv \circ
\foata^{-1}$; hence, the resulting right-subtree of the maximal element in the heap of $\varphi
(w)$ is simply $A$.  Whereas the subtree of the original right subword $B=75826$ has become
the left-subtree of the heap of $\varphi (n)$, but with $\varphi$
applied to it.   \textcolor{black}{The action of $\rev$ switches left and
right children at each vertex, while $\foata \circ \inv \circ
\foata^{-1}$ switches back the children at every node except the
records.  In particular, if $B=CmD$ with maximum $m$, then $m$ becomes
a record of $\rev{w}$ so $\foata \circ \inv \circ
\foata^{-1}$ undoes the swapping at every vertex derived from $C$,
while the vertices derived from $D$ are restored accordingly as they
would have been were $n$ and $A$ empty.}  The net effect is that $\varphi (AnB) = \varphi (B)nA$.  

The resulting heap will have the same underlying \emph{unlabeled graph structure} as the
original heap, but does not preserve the \emph{binary tree} structure, since some left edges
become right edges, while others stay the same.    
\end{proof}
\begin{defn}\label{def:treeGroup}
Let $T$ be unlabeled binary tree with $n$ vertices. For every vertex $v$ of $T$, let $\tau_{v}$ be the map
that interchanges (toggles) the left and right subtrees of $T$.  Let $\Gamma (T) =\langle
\tau_{v}:v\in T\rangle$ be the subgroup generated by these involutions, which we can think
of as a subgroup of $\frakS_{n}$ of size dividing $2^{n}$.
\end{defn}
Note that $\varphi$ acts as an element of $\Gamma (T)$, where $T = T(H(w))$ is the
underlying unlabeled tree of the heap of a permutation.  

The recursive structure described in the above lemma allows us to prove interesting facts
about the orbit structure and that certain natural statistics are homomesic.  The reader is
invited to check each statement against the two orbits of $\overline{\varphi}$ displayed (with their
heaps) in Figure~\ref{fig:phiOrbits}. To view this as an action of $\varphi$ itself, just
drop all parentheses from the listed permutations.  

\begin{thm}\label{thm:phi}
The orbits of the action of $\overline{\varphi}$ (or $\varphi$) on $\frakS_{n}$, satisfy the
following properties.  

\begin{enumerate}

\item The size of each $\varphi$-orbit (equivalently
  $\overline{\varphi}$-orbit) is a power of 2.  Specifically if $w$
  lies in the orbit, define the \textbf{height} of the heap $H(w)$ to
  be the number of edges $h$ in a maximal path from the root (to a
  leaf); then the size of the orbit is
  $2^{h}$. \textcolor{black}{Consequently, as noted by Mike Joseph, the
    GCD of the orbit sizes is the greatest power of 2 that is less
    than or equal to $n$; for if $2^k<n$, then any binary tree
    containing $1, 2, ..., n$ has at least one vertex that is at least
    $k$ steps away from the root.}
\item Let $\Fix w$ denote the number of fixed points, i.e., 1-cycles, of $w$.  Then the
statistic $\Fix$ is 1-mesic with respect to the action of $\overline{\varphi}$. 
Equivalently, $\Rasc = $\#record-ascents (see below) is 1-mesic with respect to the action of $\varphi$.

\item For fixed values $i \neq j$ in $[n]$, let $\II_{i<j}(u)$ denote the indicator statistic
of whether $i$ occurs to the left of $j$ in the one-line notation of $u$.  Then $\II_{i<j}$ is
$\frac{1}{2}$-mesic with respect to the action of $\varphi$.  

\item Similarly for fixed $i\in [n]$, let $\II_{(i, n)}$ denote the indicator statistic of
whether $i$ and $n$ lie in the same cycle of $w$.  Then $\II_{(i,n)}$ is $\frac{1}{2}$-mesic
with respect to the action of $\overline{\varphi}$.  
\end{enumerate}
\end{thm}
\begin{proof}
The proof of each statement proceeds by induction, using the lemma above, and assuming the
statement is true for the action of $\varphi$ on $\frakS_{k}$ for every $k<n$.  

%Suppose first that $w=AnB$ where neither $A$ nor $B$ is empty. 

1. If $n=1$, then there is only one orbit of size 1 (the base case).  Otherwise, $w=AnB$
where at least one of $A$ or $B$ is nonempty.  
By Lemma~\ref{lem:recurPhi}, the action of $\varphi$ on $w$ looks as follows: 
\begin{equation}\label{eq:orbPhi}
AnB \mapsto \varphi (B)nA \mapsto \varphi (A)n\varphi (B) \mapsto \varphi^{2}(B)n\varphi (A)
\mapsto \varphi^{2}(A)n\varphi^{2}(B)\mapsto \dotsb 
\end{equation}
In order for the action to return to $AnB$, we must apply the action an even number of
times, so the correct set of values lies to the left of $n$.  The first time these values will be
the same partial permutation as $A$ is the size of the orbit of $A$ in $\frakS_{\#A}$, which by inductive
hypothesis is a power of 2, say $2^{a}$.  (If $A$ is empty, take $a=0$.)  Similarly for the
values in $B$, which will return to their initial order after say $2^{b}$ steps.  Thus, the
size of the orbit of $w=AnB$ is $2\cdot \lcm (2^{a}, 2^{b}) = 2^{\max \{a,b\}+1}$.  The heap
of $w$ has root $n$ and subtrees $H(A)$ and $H(B)$.  
So assuming inductively that $a$ represents the height of the heap of $A$ and $b$ that of $B$, 
the height of the heap of $w$ is $\max \{a,b \} + 1$.  

Note that in the
case where exactly one of $A$ or $B$ is empty, we get that the orbit of $AnB$ is exactly
twice as long as that of the nonempty subword, e.g., 
\[
An \stackrel{\varphi}{\mapsto}
nA \stackrel{\varphi}{\mapsto}
\varphi (A)n \stackrel{\varphi}{\mapsto}
n \varphi (A) \stackrel{\varphi}{\mapsto}
\varphi^{2}(A)n \stackrel{\varphi}{\mapsto} 
n\varphi^{2}(A)  \stackrel{\varphi}{\mapsto} \dotsb . 
\]

2. First note that the statement holds trivially for $n=1$ and $n=2$, where there is only a
single orbit of $\varphi$.  Also note that the statistic $\Fix w$ translates to $\Rasc
\foata (w)$ where for $u\in \frakS_{n}$ given in one-line notation
\[
\Rasc u : = \#\{\text{records of } u \text{ immediately followed by another record or in the
final position} \}.
\]
In other words, $\Rasc u$ counts the number of records that are also \emph{ascents}. 
For example, for  $w = 847296513 = (4 2)(6)(8 1)(9 3 7 5)$ as in 
Equation~(\ref{eq:foata}), $\Fix w = 1$, and the 1-cycle $(6)$ translates
into the only record of $\foata (w) = 426819375$ which is immediately followed by a larger entry.  
In particular, a nonempty substring to the right of $n$ contributes nothing to $\Rasc u$
(equivalently in $\foata^{-1}(u)$ there are no 1-cycles to the right of $n$.)

We want to show that $\Rasc$ is 1-mesic on orbits of $\varphi$.  
Arguing as above, we claim that for $w=AnB$
\begin{align}\label{eq:Rasc}
\begin{split}
\text{Average of }\Rasc \text{across a }\varphi \text{-orbit}  = 
&\phantom{+}\frac{1}{2}\cdot \text{Average of }\Rasc \text{across a }\varphi \text{-orbit} \text{ of }A\\ 
&+\frac{1}{2}\cdot \text{Average of }\Rasc \text{across a }\varphi \text{-orbit} \text{ of }B\\  
&+\frac{1}{2} \text{, \textbf{in the case that} $A$ or $B$ is empty}.  
\end{split}
\end{align}
For the $\varphi$-orbit of $w$ simply alternates the orbits of the substrings $A$ and $B$ to
the left of $n$, with substrings to the right of $n$ making no contribution.  So each orbit
contributes half of its average to the total.  If one of the substrings is empty, then half
the time $n$ will contribute 1 to $\Rasc$, otherwise $n$ never contributes.  By induction
hypothesis, the average of $\Rasc$ across any orbit (hence superorbit) for either substring
is one, and the result follows from (\ref{eq:Rasc}).  

3. The statement is clear for the single $\varphi$-orbit when $n=2$.  WLOG assume $i<j$.  
If $i,j$ lie in different substrings, or if
$j=n$, then clearly $i$ alternates between being before or after $n$ by
Lemma~\ref{lem:recurPhi}.  The remaining case is that $i$ and $j$ live in the same
substring, say $A$, so their relative order is determined by their relative order in the
substrings with the same elements as $A$ in Equation~(\ref{eq:orbPhi}), i.e., $A, A, \varphi (A),
\varphi (A), \varphi^{2}(A), \varphi^{2} (A)\dots $.  By induction, $\II_{i<j}$ is
$\frac{1}{2}$-mesic for the action of $\varphi$ on $A$, so the same is true for the ``orbit''
where each element is repeated twice; equivalently, $\II_{i<j}$ is $\frac{1}{2}$-mesic for
the $\varphi$-orbit of $w=AnB$.  

4.  Clearly $i$ alternates between being before or after $n$ in the one line notation
of $\foata u = w = AnB$ by Lemma~\ref{lem:recurPhi}.  This translates into an alternation of
$i$ and $n$ being in the same cycle in $u$, since every element to the right of $n$ in $w$
ends up in the cycle with largest element $n$ when $\foata^{-1}$ is applied.  
\end{proof}

% \begin{cor}\label{cor:phiFoatic}
% All the above statements also hold for the Foatic map $\overline{\varphi}$ acting on
% $\frakS_{n}$.  
% \end{cor}

Data on orbit sizes for $\varphi$ are given in Table~\ref{tab:Drev-inv}.  

\begin{table}
\caption{Data on orbit sizes for reversal-inversion}
\label{tab:Drev-inv}
\bigskip

\begin{tabular}{l|c|c|c|c|c|c|c|c|c|c|c|}
     \phantom{CENTER}$n$ & 1  & 2 & 3 & 4 &  5 &   6  &  7 &     8&     9 &     10 &      11    \\ \hline 
    \# of orbits:        & 1  & 1 & 2 & 5 & 19 &  84 & 448 &  2884& 21196 & 174160 & 1598576    \\   
  LCM of orbit sizes:    & 1  & 2 & 4 & 8 & 16 &  32 &  64 &   128&   256 &   512  &  1024      \\ 
  GCD of orbit sizes:    & 1  & 2 & 2 & 4 &  4 &   4 &   4 &     8&     8 &     8  &     8      \\ 
  Longest orbit size:    & 1  & 2 & 4 & 8 & 16 &  32 &  64 &   128&   256 &   512  &  1024      \\ 
 Shortest orbit size:    & 1  & 2 & 2 & 4 &  4 &   4 &   4 &     8&     8 &     8  &     8      \\ 
 Size of $\id$'s orbit:  & 1  & 2 & 4 & 8 & 16 &  32 &  64 &   128&   256 &   512  &  1024      \\ 
\hline
\end{tabular}
\end{table}

\section{Complement-inversion}\label{sec:comp-inv}

The next nicest Foatic action appears to be the one given by: 
\begin{equation}\label{eq:gamma}
\overline{\gamma}: \frakS_{n} 
\stackrel{\foata}{\rightarrow} 
\frakS_{n} 
\stackrel{\comp }{\rightarrow}
\frakS_{n} 
\stackrel{\foata^{-1}}{\rightarrow}
\frakS_{n} 
\stackrel{\inv}{\rightarrow} 
\frakS_{n}. 
\end{equation}
%This one has orbit sizes which are all even and conjecturally divide $2n!$ (perhaps even
%$2(n-1)!$). The fixed point statistic is also homomesic for $n\leq 11$.

\begin{eg}\label{eg:gamma} Let $w=(2)(3)(514)\in \frakS_{5}$ in CCD.  Then the successive action of
$\overline{\gamma}$ on $w$ is detailed below: 
%%% Note that alignat alternates between flushright and flushleft, hence some & need to be doubled
\begin{alignat*}{5}\label{eqn:phiEG}
w = &(2)(3)(514) && \mapsto 23514 \mapsto 43152 \mapsto (431)(52)   &&\mapsto (413)(52)   && = \overline{\gamma}(w)\\ 
    & (413)(52)  && \mapsto 41352 \mapsto 25314 \mapsto (2)(5314)  &&\mapsto (2)(5413)  && = \overline{\gamma}^{2}(w)\\
    & (2)(5413) && \mapsto  25413 \mapsto 41253 \mapsto (412)(53)   &&\mapsto (421)(53)	 && = \overline{\gamma}^{3}(w)\\
    &  (421)(53) && \mapsto 42153 \mapsto 24513 \mapsto (2)(4)(513) &&\mapsto (2)(4)(531) && = \overline{\gamma}^{4}(w)\\ 
    &(2)(4)(531) && \mapsto 24531 \mapsto 42135 \mapsto (4213)(5) &&\mapsto (4312)(5) && = \overline{\gamma}^{5}(w)\\
    &  (4312)(5) && \mapsto 43125 \mapsto 23541 \mapsto (2)(3)(541) &&\mapsto (2)(3)(514) && = \overline{\gamma}^{6}(w)\\ 
\end{alignat*}
This example also displays (down the second column) the conjugate orbit of $\gamma$, also of size 6.  
\[
23514 \stackrel{\gamma}{\rightarrow}
41352 \stackrel{\gamma}{\rightarrow}
25413 \stackrel{\gamma}{\rightarrow}
42153 \stackrel{\gamma}{\rightarrow}
24531 \stackrel{\gamma}{\rightarrow}
43125 \Lsh 
\]
\textcolor{black}{There are $6$ fixed points distributed between the $6$
  permutations in this orbit:  $2$ is fixed by three permutations, and
  $3$, $4$, and $5$ are fixed once each.}
\end{eg}

\begin{table}
\caption{Data on orbit sizes for complement-inversion}
\label{tab:Dcomp-inv}
\bigskip

\scalebox{.95}{%
\begin{tabular}{l|c|c|c|c|c|c|c|c|c|c|c|}
     \phantom{CENTER}$n$          & 1  & 2 & 3 &  4 &  5&   6 &   7 &     8&     9 &       10 &        11    \\ \hline 
    \# of orbits:                 & 1  & 1 & 2 &   5& 15&   60&  288&  1656&  11028&    84042 &    717700    \\   
  LCM of orbit sizes:             & 1  & 2 & 4 &  24& 48&  480& 2880& 40320& 241920& 50803200 & 101606400    \\
  GCD of orbit sizes:             & 1  & 2 & 2 &   2&  2&    2&    2&     2&     2 &     2    &         2    \\ 
  Longest orbit size:             & 1  & 2 & 4 &   8& 16&   32&   80&   144&   360 &  1260    &      2880    \\ 
 Shortest orbit size:             & 1  & 2 & 2 &   2&  2&    2&    2&     2&     2 &     2    &         2    \\ 
 Size of $\id$'s orbit:           & 1  & 2 & 4 &   8& 16&   32&   64&   128&   256 &   512    &      1024    \\ 
\hline
\end{tabular}%
}
\end{table}

\textcolor{black}{We examined all orbits of the action of $\gamma$ on
  $\frakS_n$ for $n\leq11$ and the number of fixed points was
  $1$--mesic in every orbit we encountered.  We also observed some
  interesting patterns to the orbit sizes, but were unable to prove
  any in generality.  Data on orbit sizes in tabulated in
  Table~\ref{tab:Dcomp-inv}.}

\textcolor{black}{The prime factorizations of the LCM of orbit sizes
  suggest a simple formula and hint at an underlying structure to the
  action, but we have been unable to make the pattern precise or
  identify such structure.
  \begin{align*}
    1&= 1  &   24&=2^3\cdot3 &     2880&=2^6\cdot3^2\cdot5 &     50803200&=2^9\cdot3^4\cdot5^2\cdot7 \\
    2&=2   &   48&=2^4\cdot3 &     40320&=2^7\cdot3^2\cdot5\cdot7 &
                                                                    101606400&=2^{10}\cdot3^4\cdot5^2\cdot7 \\
    4&=2^2 &  480&=2^5\cdot3\cdot5 &     241920&=2^8\cdot3^3\cdot5\cdot7 &
  \end{align*}
  The powers of $2$ in the LCMs are achieved by the orbit containing the
  identity permutation in every case we have examined, but we do not
  have a simple way to show this is never exceeded, nor to explain the
  other prime factors.
}

\textcolor{black}{Notice that when $n\geq2$, the numbers $1$ and $n$ are
  in the same cycle of $w$ if and only if they are in different cycles
  of $\gamma(w)$.  It follows that every orbit has even length, and
  the indicator statistic $\II_{(i, n)}$ is $\frac12$-mesic.  The GCD
  of the orbits is completely accounted for by noting that
  $\rev(123\cdots{}n)$ is in the unique $2$-cycle created by the
  action of $\gamma$ for every $n\geq2$.
}

\begin{conj}\label{conj:gamma}
The action of $\overline{\gamma}$ on $\frakS_{n}$ has the following properties.  
\begin{enumerate}

\item The statistic $\Fix$, which counts the number of fixed points (1-cycles), is 1-mesic.

\end{enumerate}
\end{conj}

\section{Complement-rotation}\label{sec:comp-rot}

\textcolor{black}{Our search of orbits for small $n$ also failed to
  produce a counter-example to the conjecture that $\Fix$ is 1-mesic
  for the Foatic action given by:  }
\begin{equation}\label{eq:rho}
\overline{\rho}: \frakS_{n} 
\stackrel{\foata}{\rightarrow} 
\frakS_{n} 
\stackrel{\comp }{\rightarrow}
\frakS_{n} 
\stackrel{\foata^{-1}}{\rightarrow}
\frakS_{n} 
\stackrel{\rot }{\rightarrow} 
\frakS_{n}. 
\end{equation}
\textcolor{black}{Our observations about orbit sizes for $\rho$ are
  tabulated in Table~\ref{tab:Dcomp-rot}.  A consequence of our choice
  of computing environment means that this tabulation is incomplete,
  since the LCM of orbit sizes could not be carried out using
  PostScripts native integers, but the prime factorization in the case
  $n=7$ is sufficient to exclude any simple structures of the form
  discussed in the preceding two section.}

\begin{table}
\caption{Data on orbit sizes for complement-rotation}
\label{tab:Dcomp-rot}
\bigskip

\scalebox{.92}{%
\begin{tabular}{l|c|c|c|c|c|c|c|c|c|c|c|}
     \phantom{CENTER}$n$          & 1  & 2 & 3 &  4 &   5 &     6  &     7   &     8       &     9     &       10 &       11   \\ \hline 
  \# of orbits:                   & 1  & 1 & 1 &  4 &   8 &    30	&    70	  &   300	&   716	    &  3360    &  7012	    \\   
  LCM of orbit sizes:             & 1  & 2 & 6 & 24 & 480 & 347760 &
                                                                     $\sim
                                                                     3\times
                                                                     10^{14}$
                                                                             &
                                                                               n/a & n/a & n/a & n/a   \\ 
  GCD of orbit sizes:             & 1  & 2 & 6 &  2 &   2 &     2	&     2	  &     2	&     2	    &     2    &     2	    \\ 
  Longest orbit size:             & 1  & 2 & 6 &  8 &  32 &   108	&   576	  &  2694	& 16864	    & 116168   & 1676162    \\ 
 Shortest orbit size:             & 1  & 2 & 6 &  4 &   6 &     6	&    10	  &    10	&    12	    &    12    &    14	    \\ 
 Size of $\id$'s orbit:           & 1  & 2 & 6 &  8 &  10 &    12 	&    14   &    16    	&    18     &    20    &    22      \\ 
\hline
\end{tabular}}

\medskip {\small\textit{\textcolor{black}{The LCM for $n=7$ has the exact value
    $295162920561600=(2)^6(3)^4(5)^2(7)(13)(23)(67)(109)(149)$.
    Entries of n/a reflect the fact that our initial computations were
    carried out using the PostScript programming language, and
    intermediate calculations exceeded the maximum value of an integer
    in our interpreter.  The values should be easily accessible to any
    serious reimplementation of the computation.}} }
\end{table}

\textcolor{black}{We were surprised to observe from our data that every
  orbit of $\overline\rho$ had a representative permutation in which
  $1$ appeared as a fixed point, but given the number of orbits we
  examined, we conjecture that this is a general property of orbits of $\overline\rho$.}

\begin{conj}\label{conj:rho}
The action of $\overline{\rho}$ on $\frakS_{n}$ has the following properties.  
\begin{enumerate}

\item The statistic $\Fix$, which counts the number of fixed points (1-cycles), is 1-mesic.  

\item Each orbit contains at least one permutation with 1 as a fixed point (aka $(1)$ as a
  1-cycle).

\item \textcolor{black}{There is a natural indicator statistic that is $\frac12$-mesic
    and accounts for the orbits having even size.}

\end{enumerate}
\end{conj}

\section{Reversal-rotation}\label{sec:rev-rot}

\textcolor{black}{The only remaining Foatic action for which $\Fix$ is
  potentially 1-mesic is}
\begin{equation}\label{eq:tau}
\overline{\tau}: \frakS_{n} 
\stackrel{\foata}{\rightarrow} 
\frakS_{n} 
\stackrel{\rev}{\rightarrow}
\frakS_{n} 
\stackrel{\foata^{-1}}{\rightarrow}
\frakS_{n} 
\stackrel{\rot}{\rightarrow} 
\frakS_{n}. 
\end{equation}
\textcolor{black}{Our data on the sizes of orbits is tabulated in Table~\ref{tab:Drev-rot}.}

\begin{table}
  \caption{Data on orbit sizes for complement-rotation}
\label{tab:Drev-rot}
\bigskip
\scalebox{.82}{%
\begin{tabular}{l|c|c|c|c|c|c|c|c|c|c|c|}
  \phantom{CENTER}$n$ & 1 & 2 & 3 & 4 & 5 & 6 & 7 & 8 & 9 & 10 & 11 \\
  \hline
  \# of orbits: &1 & 1 & 1 & 4 & 8 & 26 & 50 & 222 & 378 & 2356 & 3634 \\
  LCM of orbit sizes:& 1 & 2 & 6 & 24 & 4680 & 155195040 & 1768492025501160960 & n/a & n/a & n/a & n/a \\
  GCD of orbit sizes:&1 & 2 & 6 & 2 & 2 & 2 & 2 & 2 & 2 & 2 & 2 \\
  Longest orbit size:&1 & 2 & 6 & 8 & 36 & 102 & 726 & 1216 & 9656 & 20050 & 160128 \\
  Shortest orbit size:&1 & 2 & 6 & 4 & 6 & 4 & 8 & 8 & 12 & 12 & 12\\
  Side if $\id$'s orbit:&1 & 2 & 6 & 6 & 36 & 28 & 70 & 96 & 4312 &
                                                                    3784 & 19784 \\
  \hline
\end{tabular}%
}
\end{table}

\begin{conj}\label{conj:tau}
  The statistic $\Fix$, which counts the number of fixed points
  (1-cycles), is 1-mesic for the action of $\tau$ on $\frakS_{n}$.
\end{conj}

\section{Other Foatic maps }\label{sec:otherMaps}

\textcolor{black}{The remaining 21 Foatic actions as described in
  Section~\ref{ss:dyn} all fail to have $\Fix$ as a $1$-mesic
  statistic and fall outside of the scope of our initial
  investigation.  The dynamics of these action have been further
  studied, however, by Elizabeth Sheridan Rossi in her doctoral
  dissertation~\cite{esrthesis}, where she considered a wider range of
  statistics, as outlined in Section~\ref{ss:summ}.}

\section{Data Availability}

Part of this investigation was accomplished using programs written in
the PostScript programming language.  PostScript is not usually
considered suitable for serious mathematical research, but it was
sufficiently flexible for the early stages of this project.  One such
program is
currently available from the first author's webspace\\
\url{https://math.mit.edu/~malacroi/permHomomesyFoaticFP/FoaticActions.txt}\\
It is designed to be run using the \emph{GhostScript} interpreter, and
could be
invoked as follows.\\
\verb!gs -dNODISPLAY -dusecycles=true -dmaxpl=8 FoaticActions.txt!\\
\verb!gs -dNODISPLAY -dshoworbits=true -dmaxpl=6 FoaticActions.txt!\\
\verb!gs -dNODISPLAY -donlygood=true -dmaxpl=10 FoaticActions.txt!

In our initial investigation, we permitted $\rot[]$ and $\rot[3]$ to
take the roles of $\calA$ and $\calB$, so our data involves 49
distinct actions.
Some output is also available for interested parties.  A tabulation of
all of the orbits acting on permuations of order at most 6 can be found at\\
\url{https://math.mit.edu/~malacroi/permHomomesyFoaticFP/AllOrbits1-6.txt}\\
While a corresponding tabulation of the 4 actions conjectured to have
$\Fix$ as a homomesy can be found at\\
\url{https://math.mit.edu/~malacroi/permHomomesyFoaticFP/GoodOrbits1-8.txt}

%\section{AUTHOR NOTES (NOT FOR PUBLICATION!)}

%\subsection{Gleanings}\label{ss:glean}

%%%Also reference paper by Tewodros Amdeberhan (tamdeber@tulane.edu), which studies iterates
%%%of the other one.  This is more important for Elizabeth's thesis.  

%\subsection{ToDo}\label{ss:todo}
%
%
%\begin{enumerate}
%
%\item Flesh out the various sections as written. 
%
%\item Should we try to create a link or somehow archive the data file for interested
%parties?  
%
%\item What's the right conjecture for the very round orbit sizes for complement-inversion?  
%
%\item Add a framework for the nice Foatic actions.
%\end{enumerate}

\bibliographystyle{amsplain-ac}
\bibliography{bibliography-foatic}

%%%%%%%%%%%% End OF DOCUMENT %%%%%%%%%%%%%%%%%%%%%%%
%
\end{document}